\documentclass[11pt]{amsart}
\usepackage{geometry}                
\usepackage{graphicx}
\usepackage{amssymb}
\usepackage{epstopdf}
\usepackage{tikz}
\usepackage{multicol}
\usetikzlibrary{positioning}
\title{The legacy of Bletchley Park on UK mathematics}
\author{Daniel Shiu}

\begin{document}
\maketitle
\section{Introduction}
The second world war saw a major influx of mathematical talent into the areas of cryptanalysis and cryptography. This was particularly true at the UK's Government Codes and Cypher School (GCCS) at Bletchley Park. The success of introducing mathematical thinking into activities previously dominated by linguists is well-studied, but the reciprocal question of how the cryptologic effort affected the field of mathematics has been less investigated. Although their cryptologic achievements are not as celebrated as those of Turing, Tutte and Welchman, Bletchley Park's effort was supplemented by more eminent mathematicians, and those who would achieve eminence and provide leadership and direction for mathematical research in the United Kingdom. Amongst their number were Ian Cassels, Sandy Green, Philip Hall, Max Newman and Henry Whitehead. This paper considers how the experience of these and other mathematicians at Bletchley Park may have informed and influenced the mathematics that was produced in their post-war careers.

Three main themes will be pursued. Firstly, how the successful use of Bayesian statistics at GCCS increased the use of this perspective by post-war statisticians over the classical frequentist interpretation that was prevalent in the UK before the war. Secondly, how Bletchley Park alumni recognised the value that computation could add to mathematical study and that this contributed to the development of strong UK research areas in computational mathematics. Thirdly, how the mathematical community at Bletchley Park influenced  younger mathematicians in their choice of study and research as well as laying the foundations for collaborations within academic mathematics.
\section{The promotion of Bayesian statistics}
Our first theme has already been noted and the reader is commended to read the excellent accounts of Zabell \cite{zabell-BP},\cite{zabell-good}. Zabell notes that the Bayesian interpretation of statistics was in abeyance in the United Kingdom leading up to the war. Because statistics was relatively under-represented at Bletchley Park, mathematicians less trained in statistics, such as Good and Turing, approached the subject with fresh eyes and were able to demonstrate that the Bayesian perspective was remarkably effective when applied to cryptanalytic problems. After the war, Turing was more interested in the development of computers, but Good realised that the isnights of Bayesian statistics could have a broader application than cryptanalysis. Zabell goes on to relate how Good became a prominent advocate for the Bayesian viewpoint starting with his 1950 book \emph{Probability and the Weighing of Evidence} \cite{good1950probability}, and continuing with writings such as \emph{Rational decisions} \cite{good1952rational},\emph{The estimation of probabilities: An essay on modern Bayesian methods} \cite{good1966estimation}, \emph{Good thinking: The foundations of probability and its applications} \cite{good1983good}, and \emph{The Bayes/non-Bayes compromise: A brief review} \cite{good1992bayes}.

It is worth emphasising the importance of his Bletchley Park experiences in Good's first book.  In the preface, Good readily acknowledges his Bletchley Park colleagues (without the GCCS context): 
\begin{quote}Dr. A. M. Turing, Professor M. H. A. Newman and Mr. D. Michie were good enough to read the first draft (written in 1946) and I am most grateful for their numerous suggestions.\end{quote}
More than this though, Michie actively encouraged its publication when Good's own belief in the book was at a low ebb. As Good recalls in an interview with David Banks \cite{banks1996conversation}
\begin{quote}A second reason for moving from Manchester was that my book Probability and the Weighing of Evidence (completed in 1948) \cite{good1950probability} had been declined by the Cambridge University Press. Who the referee was I don't know, but it could have been Jeffreys, since Frank Smithies told me he saw the manuscript on Jeffreys's desk. I was discouraged and did nothing about the manuscript for several months at least. Then Donald Michie, a good friend of mine right from the Newmanry days, said, Why don't you send it in somewhere and let someone else do the work if you're not working on it? I said, Look, the manuscript's in Manchester and I'm in London. He replied, Take a train to Manchester and come back the same day on another train. I didn't have enough faith in myself without his moral support. And I did what he advised, and submitted the manuscript to Charles Griffin, who published it in 1950. (Most publication was slow at the time.) I think it wasn't such a bad little book. Somebody phoned me recently and called it a classic.\end{quote}

Although now well-established (with over 1500 citations), the initial reception for the book explains some of Good's doubts and the importance of Michie's intervention. Reviews varied from the dismissive \cite{david1951probability}
\begin{quote}It is fair to say that statisticians in general will not accept Dr Good's approach to statistical problems\end{quote}
to the qualified \cite{ahcopeland}
\begin{quote}This book is interesting reading whether or not one agrees with the author. It is a significant contribution to a field in which contributions are badly needed.\end{quote}

It is also worth noting that although the big names of British statistics at the time were not present at Bletchley Park, two other cryptanalysts became statisticians of note. Although they did not have Good's evangelical zeal for the Bayesian interpretation, we can deduce that the statistical community was seeded with other ex-Bletchley Park cryptanalysts who knew of its classified value. The first is the US statistician Solomon Kullback, known best for Kullback-Leibler divergence. Kullback was seconded to Bletchley Park for three months and was exposed to the work on Enigma and TUNNY \cite{kullbacknsa} and Good reports \cite{good1980contributions} that Kullback was a user of weighted evidence ideas. Another is E. H. Simpson of Simpson's Paradox fame. In a 2010 paper \cite{simpson} Simpson explains how in addition to its use on the TUNNY and Enigma challenges, Bayes theorem was also applied to less celebrated work on the Japanese masked JN25 codebook. The exchange of ideas between the different areas is unclear. Simpson states \cite{simpson}:
\begin{quote}Because of the compartmenting of all the various teams at Bletchley Park, I knew nothing of their [Turing, Alexander, Good] work at the time.\end{quote}
Conversely, in his paper on Turing's wartime work on statistics, in the section on repeat rates, Good notes \cite{good1979studies}
\begin{quote}E. H. Simpson and I both obtained the notion from Turing. \end{quote}

Whether or not others' knowledge of the cryptanalytic efficacy of Bayes influenced the opinion of the statistical community, it is clear that Good's mission to promote Bayesian thinking was successful. His obituary in the London Times records \cite{goodtimes}
\begin{quote}To statisticians, Good is one of the founding fathers of Bayesian statistics [...] This approach -- the Bayesian paradigm, as it is now called -- was little used before Good's work but was given an important boost by his 1950 book and his extensive subsequent writings, and is firmly established today. \end{quote}

\section{An age of experimental computational mathematics}
We now turn to the effect on mathematics of computational machinery and how Bletchley Park alumni were central to this evolution. We will not dwell long on the design and development of computers, though Bletchley Park personnel were again pivotal in this story. There were three significant computer construction projects in the immediate post-war period: the Automatic Computing Engine project, designed by Turing and supported by him 1945-47; the Manchester Computing Machine Laboratory established by Max Newman and where he was joined Jack Good and David Rees (and Turing in 1948); the Electronic Delay and Storage Automatic Calculator at Cambridge's mathematical laboratory. These projects initially produced devices (Pilot ACE, Manchester Baby, EDSAC) that are considered prototypes for later ``commercial computers'' (DEUCE, Ferranti Mk 1, Atlas). Nevertheless, they already produced computations of significance for mathematics.

The first generation machines were not used by a broad research community, and many programmes were proof-of-concept demonstrations of known mathematical results. Alongside investigations of Mersenne primes, solution of differential equations and chess programmes, the most ambitious use of the Manchester Mk 1 was probably research by Alan Turing. Prior to the war, Turing had become interested in the computations of zeroes of the Riemann $\zeta$-function \cite{turing1945method}, in particular whether the famous Riemann hypothesis could be disproved in such a way. His pre-war proposals were mechanical in nature, with a \textsterling 40 grant from the Royal Society enabling him to construct some of the gears necessary for the computation. The project took a back burner while Turing worked at Bletchley, but the Manchester Mark 1 gave him an opportunity to revisit the idea with electronic computation. The results are recorded in his 1953 paper \emph{Some calculations of the Riemann zeta-function} \cite{turing1953some}
\begin{quote}In June 1950 the Manchester University Mark 1 Electronic Computer was used to do some calculations concerned with the distribution of the zeros of the Riemann zeta-function. It was intended in fact to determine whether there are any zeros not on the critical line in certain particular intervals.\end{quote}
Although Turing's paper expresses only limited satisfaction with the results, it is a feat that few others could have accomplished and the outstanding mathematical result of first generation machinery.

Early computers were unreliable and required expert oversight. The absence of not only of ``high level'' programming languages, but even a recognisable operating system left a high bar for those wishing to engage with the new capabilities. A certain level of faith in the outcome was necessary to invest the effort, but this faith was considerably easier for those who had seen the successful use of electronic calculation on the problems of Bletchley Park. One influential voice was J. W. S. Cassels, who although he did not programme the early computers himself, encouraged their use. In his overview paper \emph{Computer-aided Serendipity} \cite{cassels1995computer}, he notes
\begin{quote}I am not a computer expert, but I am fortunate in having friends who are, and in having been a witness of early developments. The first approaches to electronic computing were made during the war with specialised devices for ballistic or cryptographic purposes\end{quote}
\begin{quote}[...]\end{quote}
\begin{quote} There was no such profession as computer science in those days. Several of my fellow research students in number theory found the challenges of the new discipline congenial and moved into it. I might add that I have always encouraged my research pupils to learn about computers.\end{quote}
One of Cassels's students, Bryan Birch, was an eager adopter of computer assistance. The famous Birch-Swinnerton-Dyer conjecture (one of the Clay Institute's Millennium Problems) motivated by computations on EDSAC2 in the sixties is a notable example. Moreover, Birch supervised a very large number of doctoral students and laid the foundations for a national strength in computational number theory that endures to this day in both academic and classified communities. Birch was an early researcher at the Atlas Computer Lab in Harwell near Oxford. Here mathematics researchers were able to perform experiments using a second generation computer, taking advantage of both hardware developments in the form of transistors as well as usability features such as the Atlas Supervisor (arguably the first modern operating system) and high level language for programming.

Unsurprisingly, the expertise of Bletchley Park alumni was attractive to the Atlas Lab and Jack Good and Donald Michie were recruited in short order. Good recalls \cite{banks1996conversation}
\begin{quote}
For a few months, the Atlas was the fastest computer in the world, but IBM overtook it. The laboratory was sixteen miles away from Oxford, along a road dangerous for an inexperienced driver, so I didn't visit the lab as often as I would have liked. I did most of my work at the college, some of which was on the ``underware'' for a ``five-year plan'' for chess programming. One day I came into the laboratory and found my office taken over...
\end{quote}

The most impressive story of Atlas and a Bletchley Park alumnus is the work of A. O. L. Atkin. In his \emph{Tribute to Oliver Atkin} \cite{birchonatkin}, Birch is unequivocal in drawing the connection between Atkin's work at then Newmanry at Bletchley Park and his subsequent successes.
\begin{quote}
       So by the time he was 20, Oliver had had more experience of a group of extremely intelligent people working together than most people ever see in their lifetimes; one can hardly imagine a more exciting mathematical environment. Oliver had learnt what mathematics was about, and could see how important computation could become.
\end{quote}       
\begin{quote}[...]\end{quote}
\begin{quote}
He wanted to know what these congruences looked like, and after Bletchley Park he thought he could use computers to find out. His opportunity occurred when the ATLAS laboratory advertised a research fellowship in 1961. [I should explain that the ATLAS computer laboratory was attached to the Atomic Energy Research Establishment at Chilton (Harwell) about twelve miles south of Oxford.] Oliver applied, and of course he was appointed; his duties were to use one of the world's most advanced computers to do research in topics of his choice. At ATLAS, he was in his element, though there were those who didn't like it when he walked down the corridor singing Wagner fortissimo       
\end{quote}

The early number theoretic successes of work at Atlas and the approach of the researchers is documented in the proceedings of 1969 Atlas Symposium \emph{Computers in Number Theory} \cite{atkin1971computers}. In their preface Atkin and Birch note
\begin{quote}The papers illustrate all aspects of the use of computers in number theory: as an essential part of a proof, as an aid to discovery [...], and negatively as a possible ally in doing what has not yet been done. The attitude sometimes maintained a few years ago, that a computer is a disreputable device in the context of ``pure'' mathematics, was noticeably absent at the symposium.\end{quote}
As well as papers by Atkin and Birch, Jack Good provided two papers to the symposium. One on congruential pseudo-random number generators has a decidedly cryptographic feel, the other was a computational investigation of the Mertens conjecture. This is a number theoretic conjecture even more powerful than the Riemann hypothesis, that by the time of the conference was already felt to be false by the number theory community. Good and his former GCHQ colleague Bob Churchhouse gather data to formulate a statistical argument for the falsity of the conjecture whilst simultaneously making a heuristic case for the Riemann hypothesis. In \cite{good1968riemann}, inkeeping with the above remarks of Atkin and Birch, they write
\begin{quote}Thus Conjecture A is not merely to be expected a priori, by mathematical common sense, but it is well supported by the numerical data. As we said before, we believe therefore that there is a good ``reason'' for believing the Riemann hypothesis, apart from the calculation of the first 2,000,000 zeros.\end{quote}
The computational perspective on Mertens conjecture proved powerful and Odlyzko and te Riele were able to produce a computer-aided counterexample in 1985 \cite{odlyzko1985disproof}.

Atkin's time at Atlas only lasted until 1970 and the lab itself only until 1975. Atkin continued to be a prolific computer and correspondent, though sparse in his writing for publication. It is worth emphasising how important some of his contributions were in the 1980s when elliptic curves began to play a major role in cryptography. The Schoof-Elkies-Atkin point counting method and Atkin-Morain primality certification are both important methods for public key cryptography. The results further underscore the symbiotic development of mathematics and cryptography. The circle is completed by the observation that the original ``non-secret encryption'' result by Cliff Cocks was informed by a seminar by Cocks's supervisor, Birch, on new computational methods for factorisation. This ushered in a cryptographic era of public key cryptography to which Atkin contributed strongly.

\section{Persistent connections}
In this last section, we consider how the working and social connections established at Bletchley Park continued as academic mathematical connections. Many of the Bletchley Park mathematicians had already known each other, indeed these pre-existing connections were often the reason for their recruitment. Welchman \cite{welchman1982hut} in particular reached recruited colleagues, students under him, those who had shared his undergraduate days and even his old school. Not all of the mathematicians knew each other however and the informal style of Bletchley Park fostered the exchange of ideas and friendships between undergraduates and established academics.

The evidence of this extension is manifested firstly in Bletchley Park alumni who completed PhDs in mathematics after the war under the supervision of other alumni. The PhD was not necessarily a pre-requisite for academic study prior to the war (for example although David Rees studied for a PhD under Welchman and Hall, he never completed the degree). After the war, the qualification became more emphasised in the UK, bringing it more into line with America and Europe. As such there was less formal process around pairing students with supervisors, and if the pair already knew each of other the study fell into place more naturally. We give some examples.

The cryptanalytic successes of Bill Tutte on TUNNY are now well-known. Mathematically he is also famous for bringing about a renaissance in the study of graph theory. It is worth noting that prior to Bletchley Park, Tutte had studied and graduated as a chemist (though one with a strong talent for solving mathematical problems), and only transferred to a study of mathematics in 1940. Very soon after, he joined Bletchley Park where he worked on the FISH ciphers, including being a member both of the Newmanry and one of the five initial members of the research section. After the war Bill Tutte was elected a Fellow of Trinity College Cambridge, in recognition of his work at Bletchley Park. He noted \cite{younger2012william}
\begin{quote}It seemed to me that the election might be criticized as a breach of security, but no harm came of it.\end{quote}
He set about obtaining a PhD in mathematics and his choice of supervisor was his colleague from the Newmanry and research section, Shaun Wylie. He had also fixed on graph theory as the topic. This was unusual as it was far from Wylie's speciality of homotopy. Supposedly a condition of the supervision was that Tutte provide his own problems. Tutte adds \cite{younger2012william}
\begin{quote}in those days, graph theory had a low reputation in the mathematical world. so it is not surprising that Shaun Wylie, as my Ph.D. supervisor, advised me to drop graph theory and take up something respectable, such as differential equations. If one assumes that graph theory was my m\'etier, it was just as well that I had the prestige of a Fellow of Trinity.\end{quote}
Tutte was prolific as a PhD student, quickly producing four papers each of which was hugely influential, with the Tutte 1-factor theorem and Tutte polynomial rapidly becoming eponymous. His thesis went further still, establishing the subject now known as matroid theory.

Whereas Tutte had very clear ideas of the problems that he wished to tackle even though they were outside the sphere of his friend and supervisor, Peter Hilton's path was less planned. He relates how he came to study topology \cite{steen}
\begin{quote}
two of my colleagues [at Bletchley Park] were Henry Whitehead and Max Newman. Due to the peculiar circumstances of the war, I was on an equal footing with them. In fact, I was simply a young man who had taken a wartime degree and they were eminent mathematicians. I became, in particular, very, very friendly with Henry Whitehead, on first-name, beer-drinking terms. After the war, Henry Whitehead invited me to come back to Oxford and be his first post-war research student. I said, ``But I don't know anything about topology." And he said, ``Oh, don't worry, Peter, you'll like it." So in fact, I didn't even know what the subject was. I went back to Oxford, and I studied topology. Whatever Henry Whitehead had specialized in, I would have studied. His personality was so attractive, that it was clear that it was going to be great fun to work with him. It turned out to be not only fun but extremely exacting and demanding - it was a marvellous experience. I really took up topology because it was Whitehead's field.
\end{quote}
Whitehead's prediction was absolutely correct. Hilton went on to a very successful research in topology as well as authoring multiple authoritative textbooks in this and other areas.

Another of Whitehead's students, Gordon Preston was even younger than Hilton (Hilton joined Bletchley Park in 1942 at the age of 18; Preston in 1944 at  the age of 19). Whereas Hilton's four terms at Oxford allowed him to claim a wartime BA, Preston undertook another two years of undergraduate study. These included popular lectures given by Whitehead and Preston went on to become his graduate student \cite{preston}, 
\begin{quote}For the first year I worked with him on algebraic topology. He then went off to the States for a year and my new supervisor was E C (Edward) Thompson and, with him, I changed my subject to algebraic geometry, with a strong emphasis on commutative ring theory and with not so much geometry.\end{quote}
Preston's mathematical interests were influenced work of other Bletchley Park colleagues\cite{preston}
\begin{quote}
My semigroup influence came because of the friends I had met at Bletchley Park. I used to go to the National Physical Laboratory at Teddington - this was while I was an undergraduate again from 1946 to 1948 - to keep in touch with the development of the pilot ACE (Automatic Computing Engine) that Alan Turing had been drafted there to develop. I was interested in the mathematical papers of my other Bletchley friends and, in particular, this was how I came to read David Rees's papers. I have no records, and my memory may be playing me false, but I believe the first paper I read on semigroups was his paper `On semi-groups' from the Proceeding of the Cambridge Philosophical Society, 1940, together with the small technical note that followed it, ``Note on semi-groups", ibid., 1941.
\end{quote}
He also notes that a significant informal resource for semigroups were lectures by Philip Hall at Cambridge; again Bletchley Park colleagues are involved \cite{preston}
\begin{quote}My introduction to them came when Sandy Green, probably in 1951, lent me a copy of his notes on the lectures given one year by Philip Hall. I made a hand-written copy of them - this was before the days of the ubiquitous photocopier - which I still have. This was exciting and tremendous stuff. Semigroups were used quite naturally in Hall's lectures, at any rate for that year.\end{quote}
Whitehead had returned to the UK when Preston submitted his thesis and was internal examiner and David Rees as external examiner \cite{preston}
\begin{quote}
My D. Phil. examiners were David Rees and Whitehead. David gave me a list of detailed comments which Henry Whitehead would not let him ask me about in the oral examination - Henry was the chairman of examiners - because Henry had to rush off to captain a cricket team playing that afternoon. Henry asked me a number of questions himself, including some about my chapter on semigroups; and in my answers to his questions which were critical questions I was able apparently to show that he had failed to grasp properly the concepts involved. After this I was invited to come and watch the cricket game, which I took as a sign that I had passed my viva.
\end{quote}

The area of semigroups in UK mathematics is dominated by Bletchley Park mathematicians. At its heart we find the eminent mathematician Philip Hall. In contrast to Whitehead's gregarious nature, Hall was more reserved\begin{footnote}{``When [the mathematician Olga] Taussky accused Hall of being the worst recluse in Cambridge, he replied `No, Turing is worse'."\cite{hallbio}}\end{footnote} However, accounts of his post-war lecturing show that he presented a compelling mathematical vision \cite{hallbio}
\begin{quote}During the next five years he gave lectures in which he explored not groups so much as aspects of all the other algebras. He was full of vigour and brimming with ideas. When he was not presenting his own results, he illuminated others with his own interpretation, always going further than they had gone. Students felt that they were hearing today what the rest of the world would only hear tomorrow.'\end{quote}
Two junior Bletchley Park colleagues ``Sandy'' Green (who arrived at Bletchley Park in 1944 at the age of 18) and Derek Taunt (who started in 1941 at the age of 23) went on to study under Hall. Green notes \cite{greenbio}
\begin{quote}At Cambridge I was much impressed by lectures given by the algebraist Philip Hall; and in a different way by D E Littlewood.\end{quote}
Taunt's college obituary notes \cite{tauntobit}
\begin{quote}
Derek remained at Bletchley until VE day, when suddenly there was nothing left to do. He was posted to Teddington to work on aerodynamics, but once the Japanese war was over he came back to Cambridge to pick up where he had left off. Hardy had retired and Derek chose instead to work with the distinguished group theorist Philip Hall, who became a life-long friend.
\end{quote}

We summarise the supervisor-student relationships of this section in the following diagram.
\bigskip

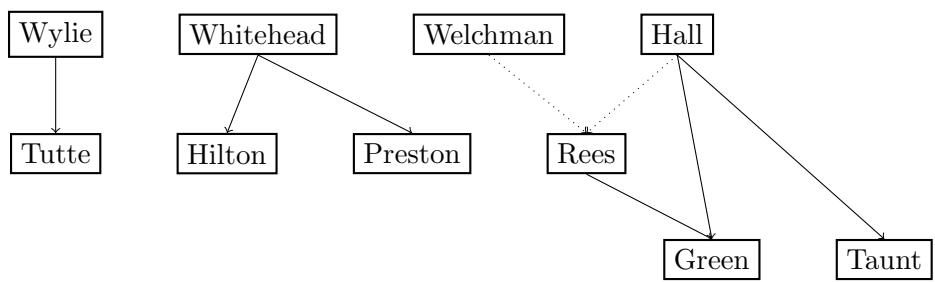
\begin{figure}[h]\centering
\begin{tikzpicture}[namenode/.style={rectangle,draw=black, thick,minimum size =5mm},]
\node[namenode]  (wylie){Wylie};
\node[namenode]  (tutte) [below=of wylie] {Tutte};
\node[namenode] (whitehead) [right=of wylie]{Whitehead};
\node[namenode] (hilton) [right=of tutte]{Hilton};
\node[namenode] (preston) [right=of hilton]{Preston};
\node[namenode] (welchman) [right=of whitehead]{Welchman};
\node[namenode] (hall) [right=of welchman]{Hall};
\node[namenode] (rees) [right=of preston]{Rees};
\node (null) [right=of rees]{};
\node[namenode] (green) [below=of null] {Green};
\node[namenode] (taunt)[right=of green]{Taunt};

\draw[->] (wylie.south)--(tutte.north);
\draw[->] (whitehead.south)--(hilton.north);
\draw[->] (whitehead.south)--(preston.north);
\draw[dotted,->] (welchman.south)--(rees.north);
\draw[dotted,->] (hall.south)--(rees.north);
\draw[->] (hall.south)--(green.north);
\draw[->] (rees.south)--(green.north);
\draw[->] (hall.south)--(taunt.north);
\end{tikzpicture}
\caption{Bletchley Park supervision graph}
\end{figure}

\bigskip

Another measure of collaboration is co-authorship of research and textbooks. Mathematics papers tend more to single authorship than other scientific areas and this was even more true before the modern era of cheap travel and rapid communications. There are however some notable joint works by Bletchley Park alumni. These papers should be judged both substantial and influential with Google scholar recording dozens of citations for each. The most cited by far (798 citations) is the definitive introductory homology textbook \emph{Homology theory: An introduction to algebraic topology} \cite{hilton1967homology}, written by Peter Hilton and Shaun Wylie. Glowingly reviewed on its publication, it became the standard text for final year undergraduate courses in topology and early graduate students. It is now regarded as a classic. 

Another work by Hilton, this time with his supervisor Whitehead \emph{Note on the Whitehead product} \cite{hilton1953note} was accepted to the Annals of Mathematics (regularly ranked as the leading mathematics journal in the world).

Hilton also co-authored a paper with David Rees \emph{Natural maps of extension functors and a theorem of RG Swan} \cite{hilton1961natural}.

There was also a significant collaboration on a family of semi-groups by David Rees and Sandy Green \emph{On semi-groups in which $xr= x$} \cite{green1952semi}.

The other co-authorship of note is between Max Newman and Alan Turing. Prior to their work at Bletchley Park, Newman and Turing had already collaborated on computational models in their 1942 paper \emph{A formal theorem in Church's theory of types} \cite{newman1942formal}. Subsequent to their war work we have noted that Newman drew Turing to Manchester to work on the Manchester baby. Public interest in the Manchester work was piqued by talk of artificial intelligence. As part of stating the case more clearly Newman and Turing took part in a BBC radio broadcast \emph{Can automatic calculating machines be said to think?} \cite{turing2004can} the discussion was moderated by Richard Braithwaite (a professor of philosophy) and an opposing view was put by Geoffrey Jefferson (a professor of neurosurgery). The broadcast was transcribed and published in \emph{Machine Intelligence} (OUP) in 1999.

We summarise the academic joint publications of Bletchley Park alumni in diagrammatic form

\begin{figure}[h]\centering
\begin{tikzpicture}[namenode/.style={rectangle,draw=black, thick,minimum size =5mm},]
\node[namenode] (hilton) {Hilton};
\node[namenode] (wylie) [left=of hilton] {Wylie};
\node[namenode] (green) [right= of hilton] {Green};
\node[namenode] (rees) [above=of green] {Rees};
\node[namenode] (whitehead) [above=of hilton] {Whitehead};
\node[namenode] (turing) [left=of wylie] {Turing};
\node[namenode] (newman) [above=of turing] {Newman};

\draw (newman.south)--(turing.north);
\draw (wylie.east)--(hilton.west);
\draw (hilton.north)--(whitehead.south);
\draw (hilton.east)--(rees.west);
\draw (rees.south)--(green.north);
\end{tikzpicture}
\caption{Bletchley Park co-authoring graph}
\end{figure}
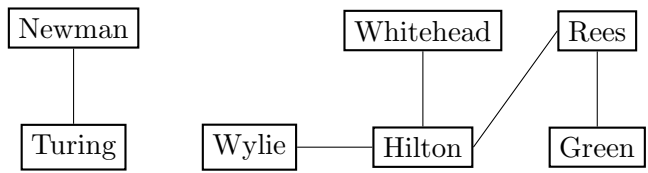

\section{Summary} We have provided evidence that Bletchley Park experiences of mathematicians went on to significantly affect their subsequent researches and publications and, as a result, shaped the direction of post-war mathematics in the United Kingdom. The effects are most prominent in Jack Good's promotion of Bayesian statistics based on his experiences with techniques used in cryptanalysis and the encouragement of his colleagues, in early computational number theory results by Turing, Good and Atkin, and in the continuation of associations begun at Bletchley Park within academia.

\bibliographystyle{plain}
\bibliography{BPmathematics}

\begin{appendix}
\section{List of mathematicians at Bletchley Park}

In the course of researching this paper, various sources have listed names of mathematicians at Bletchley Park that do not appear on the most easily accessed lists. The definition of mathematician is problematic in this case. Some of the younger mathematicians had barely begun their undergraduate careers (Bill Tutte had formally transferred to graduate study of mathematics only a few months earlier). Several people recruited due to their mathematics background remained at GCHQ and never returned to their studies. The PhD was not yet a \emph{de facto} requirement for an academic career and some notable names went on to academic careers without this qualification (David Rees for example). The following then is a list of people at Bletchley Park that can be confirmed to have enrolled in some university mathematics course before their classified work. 

\begin{multicols}{2}
\begin{itemize}
\item Hugh Alexander 
\item Michael Ashcroft \cite{preston}
\item A. O. L. (Oliver) Atkin 
\item Dennis Babbage 
\item Keith Batey  \cite{batey}
\item Howard Campaigne (US)  \cite{preston}
\item J. W. S. (Ian) Cassels 
\item Michael Chamberlain  \cite{welchman1982hut}
\item L. N. (John) Chown  \cite{chown}
\item Joan Clarke 
\item Michael Crum  \cite{preston}
\item Harold Fletcher  \cite{welchman1982hut}
\item Joe Gillis  \cite{preston}
\item I. J. (Jack) Good 
\item J. A. (Sandy) Green 
\item Philip Hall 
\item John Herivel 
\item Peter Hilton 
\item Walter Jacobs (US) \cite{good_codebreakers}
\item John Jeffreys 
\item Solomon Kullback (US) 
\item Kenneth LeCouteur \cite{good_codebreakers}
\item Tim Molien \cite{good_codebreakers}
\item Max Newman 
\item Rolf Noskwith \cite{noskwith_action}
\item Gordon Preston 
\item David Rees 
\item Margaret Rock \cite{howard2013dear}
\item Bob Roseveare  \cite{roseveare}
\item Bernard Scott \cite{banks1996conversation}
\item Edward Simpson 
\item Abraham Sinkov (US) 
\item Derek Taunt 
\item Alan Turing 
\item Bill Tutte 
\item Peter Twinn  \cite{batey2017dilly}
\item Houston Wallace  \cite{welchman1982hut}
\item Philip Watson  \cite{preston}
\item Gordon Welchmann 
\item J. H. C. (Henry) Whitehead 
\item Jimmy Whitworth  \cite{simpson}
\item Shaun Wylie 
\item Leslie Yoxall 
\end{itemize}\end{multicols}
\end{appendix}

\end{document}